\documentclass[12pt,leqno]{article}
\usepackage[francais]{babel}
\usepackage[latin1]{inputenc}
\usepackage[T1]{fontenc}
\usepackage{amsmath}
\usepackage{amsfonts}
\usepackage{amssymb}
\newcommand\Z{\mathbb{Z}}
\newcommand\C{\mathbb{C}}
\newcommand\N{\mathbb{N}}
\newcommand\R{\mathbb{R}}

\newcommand\resp{\mathrm{resp\ }}
\newcommand\rg{\rightarrow}
\newcommand\ba{\backslash}

\begin{document}

\title{Sur un ensemble de Besicovitch}
\author{Jean--Pierre kahane}
\date{}
\maketitle

Un ensemble du plan est appelé ensemble de Besicovitch s'il est d'aire nulle et contient un segment de droite dans toute direction. Si un ensemble du plan est d'aire nulle et qu'il contient un segment dans toute direction appartenant à  un certain angle, il suffit d'en prendre les  images par des rotations en nombre fini, et de réunir ces images, pour avoir un ensemble de Besicovitch. Je me permettrai dans cet article d'étendre le sens du terme, et d'appeler ensemble de Besicovitch, dans le plan ou dans l'espace $(\C,\R^2,\R^n)$ tout ensemble dont la mesure de Lebesgue  est nulle dans la dimension considérée, et qui contient des segments de droite dans toutes les directions appartenant à un ensemble ouvert de directions.

Mon point de départ sera un exemple que j'ai donné en 1959 dans  l'Enseignement mathématique, à savoir la réunion dans le plan des segments de droites dont les extrémités appartiennent à deux ensembles de Cantor de dimension $1/2$, homothétiques l'un de l'autre dans le rapport $2$, et non colinéaires \cite{Kah}. L'exemple est bon mais la justification que j'avais donnée est incorrecte. Il y a deux manières de justifier cet exemple, l'une, très rapide, par l'utilisation du théorème de projection de Besicovitch \cite{BesSam2,Matt}, et l'autre, qui entre plus profondément dans les  propriétés arithmétiques de l'ensemble, par un méthode de Richard Kenyon \cite{Keny}. Je donnerai les deux démonstrations, puis des variantes de cet exemple dans le plan et dans l'espace, et quelques commentaires.

\section{Le modèle réduit et la première démonstration}

L'exemple donné en \cite{Kah} est le transformé par affinité du modèle réduit que voici.

Dans le plan de la variable complexe $z=x+iy$, on considère les deux ensembles de Cantor
\begin{align}
E &= \left\{\displaystyle\sum_{j=1}^\infty \varepsilon_j 4^{-j}, \varepsilon_j \in \{0,1\}\right\}\\
E'&= 2E+i\,, 
\end{align}
dont les segments supports sont l'intervalle $\Big[0,\frac{1}{3}\Big]$ de la droite réelle $y=0$, et l'intervalle $\Big[0,\frac{2}{3}\Big]+i$ de la droite $y=1$. La réunion des segments de droite qui les joignent est
\begin{equation}
B=\{(1-h)z+hz'\,,\ z\in E\,,\ z'\in E'\,,\ 0\le h\le 1\}\,.
\end{equation}
Les sections horizontales de $B$ sont les
\begin{equation}
B_h=\{(1-h)z+hz'\,,\ z\in E\,,\ z'\in E'\}\ (0\le h\le 1)
\end{equation}
qu'il sera commode d'écrire $B_h =(1-h)B_u^*$ avec
\begin{equation}
B_u^* = E +\frac{u}{2} E'\,,\ \ \frac{u}{2} = \frac{h}{1-h} \quad (u>0,\ 0<h<1)\,.
\end{equation}
Les parties réelles des $B_u^*$ s'écrivent
\begin{align}
E_u &= ReB_u^* = \Big\{\displaystyle \sum_{j=1}^\infty (\varepsilon_j +u\varepsilon_j') 4^{-j}\,,\ (\varepsilon_j,\varepsilon_j')\in \{0,1\}^2\Big\}\\
&=\Big\{\sum_{j=1}^\infty a_j 4^{-j}\,,\ a_j \in \{0,1,u,1+u\}\Big\}\,.\nonumber
\end{align}

Remarquons que

\begin{align}
E_2 &= [0,1]\,,\ E_{1/2} = [0,\frac{1}{2}]\,,\nonumber\\
E_1 &=\Big\{ \sum_{j=1}^\infty \alpha_j 4^{-j}\,,\ \alpha_j \in \{0,1,2\}\Big\}\nonumber
\end{align}
donc $E_1$ est de mesure nulle et sa dimension est $\log 3/\log 4$, et
\begin{align}
E_{1/u} = \frac{1}{u} E_u\,,
\end{align}
ce qui permet de borner l'étude à $0<u<1$. Enfin
\begin{align}
4 E_u = \Big\{a_1 + \sum_{j=1}^\infty a_{j+1}4^{-j}\Big\} = \{0,1,u,1+u\} + E_u\,,
\end{align}
égalité qui va jouer un rôle essentiel dans l'étude arithmétique.

Examinons les directions des segments joignant les points $z$ de $E$ aux points $z'$ de $E'$. Ils ont tous $1$ pour hauteur, et leur projections sur l'axe réel ont pour mesures
\begin{align}
Re(z'-z) &= \displaystyle\sum_{j=1}^\infty (2\varepsilon_j' - \varepsilon_j) 4^{-j}\nonumber\\
&=\displaystyle \sum_{j=1}^\infty (2\varepsilon_j' + (1-\varepsilon_j) -1)4^{-j}\,.\nonumber
\end{align}
L'ensemble de ces mesures est $E_2 - \sum\limits_{j=1}^\infty 4^{-j}$, c'est--à--dire $\big[-\frac{1}{3}\,,\ \frac{2}{3}\big]$. Les segments ont donc toutes les directions entre $\frac{y}{x}=-3$ et $\frac{y}{x}=\frac{3}{2}$, c'est--à--dire les directions des segments joignant les intervalles--supports de $E$ et de $E'$, $\big[0,\frac{1}{3}\big]$ et $\big[i,\frac{2}{3}+i\big]$.

L'ensemble $B$ contient bien des segments de droites dans toutes les directions appartenant à un ensemble ouvert de directions.

Reste à montrer que $mes_2\, B=0$, soit $mes_1B_h=0$ pour presque tout $h$, soit $mes_1 E_u=0$ pour presque tout $u$ ($mes_\alpha$ est la mesure de Lebesgue en dimension~$\alpha$).

La preuve rapide repose sur le théorème des projections de Besicovitch. En effet, $E_u$ est la projection de $E\times E$ dans la projection qui amène $(x,y)$ en $x+uy$. Or $E\times E$ est un ensemble de dimension $1$ irrégulier selon Besicovitch, purement non--rectifiable dans la terminologie actuelle (\cite{BesSam2}, \cite {Matt} p.~204). Le théorème de Besicovitch dit justement que $mes_1 E_u=0$ pour presque tout $u$ (\cite{BesSam2},\cite {Matt} p.~250).

\section{La méthode de Kenyon, seconde démonstration}

Je vais maintenant donner une autre preuve, calquée sur une étude de Richard Kenyon \cite{Keny}. On va montrer que $mes_1 E_u=0$ quand $u$ est irrationnel ; cela suffit pour établir que $B$ est un ensemble de Besicovitch. Voici une description plus complète des propriétés des $E_u$.

Introduisons d'abord une notation. A tout entier $n$ positif, associons $n^* \in \{1,2,3\}$, le premier chiffre non nul à partir de la droite dans l'écriture de $n$ en~base~$4$ :
$$
n^* = n 4^{-j^*}  \mod 4\,,\ j^* =\sup \{j: 4^j\ \textrm{divise}\ n\}\,.
$$
Remarquons que si $\frac{p}{q}$ est une fraction irréductible, on a soit $p^*$ et $q^*$ impairs, soit $p^*+q^*$ impair.

\vskip2mm

\textsc{Théorème} 1.--- \textit{Si $u=\frac{p}{q}$ avec $p^*+q^*$ impair, $E$ contient un intervalle, et c'est l'adhérence de son intérieur. Si $u=\frac{p}{q}$ avec $p^*$ et $q^*$ impairs, $mes_1 E_u=0$ et la dimension de Hausdorff de $E_u$ est $<1$. Si $u$ est irrationnel, $mes_1 E_u=0$.} 

\vskip2mm

La démonstration reposera sur les lemmes que voici.

\vskip2mm

\textsc{Lemme} 1.--- \textit{Si $mes_1 E_u >0$, $E_u$ contient un intervalle.}

\vskip2mm

\textsc{Lemme} 2.--- \textit{Si $E_u$ contient un intervalle, c'est l'adhérence de son intérieur.}

\vskip2mm

\textsc{Lemme} 3.--- \textit{Si $E_u$ contient un intervalle, $u$ est rationnel.}

\vskip2mm

Le lemme 1 joue un rôle clé. Le lemme 2, dont la démonstration est immédiate, sera utilisé dans celle du lemme~3. La démonstration du lemme~3 fait intervenir un pavage par des translatés de $E_u$, idée ingénieuse qui était naturelle pour Kenyon.

L'implication
\begin{align}
u\ \textrm{irrationnel}\ \Longrightarrow mes_1 E_u=0
\end{align}
résulte des lemmes 1 et 3. On examinera le cas de $u$ rationnel après la démonstration du lemme~2.

\eject

\textbf{Démonstration du lemme 1}.

L'égalité (8) signifie que $4E_u$ est la réunion de quatre translatés de $E_u$, distincts si l'on suppose $0<u<1$. Comme
$$
mes_1(4E_u) \le \sum_1^4 mes_1\ \textrm{(translaté\ de\ } E_u) = 4\ mes_1 E_u\,,
$$
les intersections deux à deux de ces translatés sont de mesure $(mes_1)$ nulle ; nous dirons que ces translatés sont disjoints presque partout. En itérant, $4^nE_u$ est la réunion de $4^n$ translatés de $E_u$, qui sont deux à deux soit confondus, soit disjoints presque partout. L'ensemble des vecteurs de translations est
$$
V_n = \{0,1,u,1+u\} + 4 \{0,1,u,1+u\} + \cdots + 4^{n-1} \{0,1,u,1+u\}\,,
$$
et
\begin{align}
4^n E_u = V_n +E_u\,.
\end{align}  

Supposons $mes E_u >0$. Soit $x$ un point de densité de $E_u$. Alors la mesure de $4^n(E_u -x) \cap [-\frac{1}{2},\frac{1}{2}]$ tend vers 1 quand $n\rg \infty$, soit
\begin{align}
\lim_{n\rg \infty} mes_1\big((V_n +E_u -4^n x) \cap \big[-\frac{1}{2},\frac{1}{2}\big]\big)=1\,.
\end{align}
On peut remplacer $V_n -4^nx$ dans (11) par l'ensemble $W_n$ des $w_n \in V_n -4^nx$ tels que $w_n +E_u$ rencontre $[-\frac{1}{2},\frac{1}{2}]$, soit
$$
W_n = \big((V_n  -4^n x) \cap \big(\big[-\frac{1}{2},\frac{1}{2}\big]-E_u\big)\,,
$$
et comme les $w_n +E_u$ sont soit confondus soit disjoints presque partout, leur nombre ne dépasse pas $\dfrac{1+3\ mes_1 E_u}{mes_1 E_u}$ :
\begin{align}
\sharp\ W_n \le \frac{1+3\ mes_1 E_u}{mes_1 E_u}\, .
\end{align}
L'important dans (12) est que le second membre ne dépend pas de~$n$.

Il existe donc une sous--suite de $W_n$ qui converge, et soit $W$ sa limite. Alors (11) entraîne que
$$
mes_1 \Big((W+E_u) \cap \big[-\frac{1}{2},\frac{1}{2}\Big]\Big) = 1\,.
$$
Comme $W$ est fini, $W+E_u$ est fermé, donc
\begin{align}
W+E_u =I\,.
\end{align}
Le premier membre est une réunion finie de translatés de $E_u$, qui sont fermés. L'un au moins de ces fermés contient un intervalle (Baire). Donc $E_u$ contient un intervalle.

\vskip2mm

\textbf{Démonstration du lemme 2}

Soit $x\in E_u$, $x= \sum\limits_{j=1}^\infty
 a_j 4^{-j}$ comme en (6). Ecrivons
 $$
 x\in \sum_{j=1}^n a_j 4^{-j} + 4^{-n} E_u \,.
 $$
Si on suppose que $E_u$ a des points intérieurs, il en est de même pour $4^{-n}E_u$ quel que soit $n$, donc $x$ est limite de points intérieurs de $E_u$. Donc $E_u$ est l'adhérence de son intérieur.

\vskip2mm

\textbf{Examen du cas où $u$ est rationnel}

Les ensembles $V_n$ introduits pour avoir l'égalité (10) peuvent être définis de plusieurs façons :
\begin{align}
V_n &= \{0,1,u,1+u\} +4 \{0,1,u,1+u\} + \cdots + 4^{n-1} \{0,1,u,1+u\}\\
&=\Big\{\displaystyle\sum_{j=1}^n a_j 4^{n-j}\,,\ a_j \in \{0,1,u,1+u\}\Big\}
\nonumber\\
&=\Big\{\displaystyle\sum_{k=1}^n b_k 4^k\,,\ b_k \in \{0,1,u,1+u\}\Big\}
\nonumber\\
&=\Big\{\displaystyle\sum_{k=0}^{n-1} (\varepsilon_k+u\varepsilon_k')4^k\,,\ (\varepsilon_k,\varepsilon_k') \in \{0,1\}^2\Big\}\,.
\nonumber
\end{align}
Quand $u$ est irrationnel, il est clair qu'aucun $V_n$ n'a de point multiple. Nous allons nous intéresser à cette propriété quand $u$ est rationnel.

Commençons par une remarque générale. Le fait qu'aucun $V_n$ n'ait de point multiple peut s'exprimer ainsi : l'équation
$$
\sum \varepsilon_k 4^k + u \sum \varepsilon_k' 4^k =v \qquad (k\in \Z)
$$
où $\varepsilon=(\varepsilon_k)$ et $\varepsilon'=(\varepsilon_k')$ appartiennent à $\{0,1\}^\Z$ et sont nuls hors d'un ensemble fini, admet au plus une solution pour chaque $v\in \Z$ donné. Sous cette forme, il est évident que la propriété est invariante quand on change $u$ en $4u$ ou en~$u/4$.

Partant de la fraction irréductible $u=\frac{p}{q}$, on se ramène donc, sans modifier la propriété, au cas où $p=p^*\mod 4$ et $q=q^* \mod 4$, ce que nous supposerons maintenant. Le fait qu'aucun $V_n$ n'ait de point multiple signifie que l'équation
$$
q \sum_{k\ge 0} \varepsilon_k 4^k + p \sum_{k\ge 0} \varepsilon_k' 4^k = q \sum_{k\ge 0} \varepsilon_k'' 4^k + p \sum_{k\ge 0} \varepsilon_k''' 4^k\,,
$$
où $\varepsilon$, $\varepsilon'$, $\varepsilon''$ et $\varepsilon'''$ appartiennent à $\{0,1\}^\N$ et sont nuls pour $k$ assez grand, n'admet que les solutions $\varepsilon=\epsilon''$, $\varepsilon'=\varepsilon'''$, c'est--à--dire que l'équation
\begin{align}
q \sum_{k\ge 0} \gamma_k 4^k + p \sum_{k\ge 0} \gamma_k' 4^k =0\,,
\end{align}
où $(\gamma_k)$ et $(\gamma_k')$ appartiennent à $\{-1,0,1\}^\N$ et sont nuls pour $k$ assez grand, n'admet que la solution $(\gamma_k)=(\gamma_k')=0$. Or l'équation (15) entraîne
$$
q^* \gamma_o +p^* \gamma_o'=0 \mod 4\,.
$$
L'hypothèse que $p^* +q^*$ est impair, c'est--à--dire $p^*$ resp $q^*=2$ et $q^*$ resp $p^*=1$ ou $3$, entraîne $\gamma_o = \gamma_o'=0$, puis, en remontant à (15), $\gamma_1=\gamma_1'=0$, et finalement $\gamma_k=\gamma_k'=0$ pour tout~$k$.

\vskip2mm

\textbf{\textit{Conclusion}} : \textit{si $p^* + q^*$ est impair, aucun $V_n$ n'a de point multiple.}

\vskip2mm

Prenons maintenant comme hypothèse que $u=\frac{p}{q}$, $0<u<1$, et qu'aucun $V_n$ n'a de point multiple. Chaque $V_n$ se compose donc de $4^n$ points contenus dans l'intervalle $[0,4^n]$ et mutuellement distants d'au moins~$\frac{1}{q}$.

Soit $\mu_n$ la mesure de probabilité portée par $4^{-n}V_n$ et chargeant également tous ses points. Pour tout intervalle $J$ de longueur $|J| \ge \frac{1}{q} 4^{-n}$, on~a
\begin{align}
\mu
_n(J) \le 2q|J|\,.
\end{align}
Les $\mu_n$ convergent faiblement vers une mesure de probabilité $\mu$ portée par $E_u$ $(n\rg\infty)$, et (16) implique que $\mu$ est absolument continue, avec une dérivée $\le 2q$. Donc la transformée de Fourier de $\mu$, $\hat{\mu}$, tend vers $0$ à l'infini.~Or
\begin{align}
\hat{\mu}(t) &= \displaystyle \prod_{j=1}^\infty C(4^{-j}t)\ \textrm{avec}\ C(t) =\frac{1}{4}(1+e^{it})(1+e^{iut})\,, \\
|\hat{\mu}(t)| &= \displaystyle \prod_{j=1}^\infty |\cos 4^{-j} \frac{t}{2} |  \displaystyle \prod_{j=1}^\infty
|\cos 4^{-j} u \frac{t}{2}|\,.\nonumber
\end{align} 
Prenons $t=2q4^n\pi$ :
$$
|\hat{\mu}(t)| = \prod_{j=1}^n\! |\cos q 4^{n-j}\pi |   \prod_{j=n+1}^\infty|\!\cos q 4^{n-j} \pi  |
\prod_{j=1}^n\! |\cos p 4^{n-j}\pi |
\prod_{j=n+1}^\infty |p \cos 4^{n-j}\pi | 
\,.
$$
Les produits $\prod\limits_{j=1}^n $ valent 1, et les produits $\prod\limits_{j=n+1}^\infty$ sont indépendants de $\pi$. le dernier, qui vaut $\prod\limits_{j=1}^\infty |\cos p 4^{-j}\pi|$, n'est nul que si un facteur est nul, c'est--à--dire si $p4^{-j}\pi=\frac{\pi}{2}\mod \pi$ pour un certain $j$, c'est--à--dire si $p^*=2$. De même, le premier $\prod\limits_{j=n+1}^\infty$ n'est nul que si $q^*=2$. Donc, si $p^*$ et $q^*$ sont impairs, on~a
$$
\overline{\lim_{t\rg \infty}} |\hat{\mu}(t)|>0
$$
ce qui est incompatible avec l'hypothèse.

\vskip2mm

\textbf{\textit{Conclusion $1$ :}} \textit{si $p^*$ et $q^*$ sont impairs, l'un des $V_n$, disons $V_{n_0}$, a un point multiple, c'est--à--dire que le cardinal de $V_{v_o}$ est strictement inférieur à~$4^{n_o}$ :}
$$
\sharp V_{n_o}=\nu < 4^{n_0}\,.
$$

\vskip2mm

Poursuivons. L'égalité (10), itérée, entraîne pour $m=1,2,3,\ldots$
$$
4^{mn_o}E_u = \underbrace{V_{n_o}+V_{n_o}+\cdots +V_{n_o}}_{m\ \textrm termes} + E_u
$$
et comme $E_u \subset[0,1]$ puisqu'on a supposé $0<u<1$, l'hypothèse $u=\frac{p}{q}$ avec $p^*$ et $q^*$ impairs entraîne que, pour un $n_0$ et un $\nu<n_0$ convenables et pour tout $m$, $4^{mn_o}E_u$ est recouvert par $\nu^m$ intervalles de longueur 1, c'est--à--dire que $E_u$ est recouvert par $\nu^m$ intervalles de longueur $4^{-mn_o}$. Il en résulte
$$
\dim E_u \le \frac{\log\nu}{4\log n_o}<1\,,
$$
$\dim$ signifiant la dimension de boite (ou de Minkowski), par conséquent aussi la dimension de Hausdorff.

Revenons plus haut.

\vskip2mm

\textit{\textbf{Conclusion $2$ :}} \textit{si $p^*$ et $q^*$ est impair, $E_u$ porte une mesure de probabilité absolument continue, donc $mes_1 E_u>0$ et, d'après les lemmes~$1$ et $2$, $E_u$ est l'adhérence de son intérieur.}

\vskip2mm

La partie du théorème 1 concernant le cas $u$ rationnel est établie par les conclusions 1 et 2. Reste à démontrer le lemme~3.

\vskip2mm

\textbf{Démonstration du lemme 3}

Supposons que $E_u$ contient un intervalle ouvert, $J$. Cela impose que, quel que soit $n$, la distance de deux éléments distincts de $V_n$ soit minorée par $|J|$ ; sinon, deux translatés de $J$ par des éléments distincts de $V_n$ empiéteraient l'un sur l'autre, ce qui est impossible. Pour $N$ assez grand, $J$ contient deux translatés de $4^{-N}E_u$ disjoints
$$
v_i +4^{-N} E_u\,,\ v_i \in 4^{-N}V_N \quad (i=0\ \textrm{ou}\ 1)\,.
$$
On supposera $v_o <v_1$, c'est--à--dire
$$
J_o =v_o +4^{-N}J
$$
à gauche de
$$
J_1 =v_1 +4^{-N}J\,.
$$

On a vu que si un $V_n$ a un point multiple, $\dim E_u<1$. Donc ici tous les $V_n$ sont constitués de points distincts, à distances mutuelles $\ge |J|$. Prenant $n=2N$ dans la formule (10), on voit que $E_u$ est recouvert par $4^{2N}$ translatés de $4^{-2N}E_u$, les vecteurs translations appartenant à $4^{-2N}V_{2N}$, et l'ensemble recouvert plusieurs fois est de mesure nulle. Extrayons de ce recouvrement les recouvrements minimaux de $J$, $J_o$ et~$J_1$ :
\begin{align}
&J \subset \bigcup_{w\in W(J)} (w+4^{-2N}E_u)\nonumber\\
&J_o \subset \bigcup_{w\in W(J_o)} (w+4^{-2N}E_u)\nonumber\\
&J_1 \subset \bigcup_{w\in W(J_1)} (w+4^{-2N}E_u)\,.\nonumber
\end{align}
Observons que le recouvrement de $J_1$ est le translaté par $v_1-v_0$ du recouvrement de $J_o$. Montrons comment reconstituer le recouvrement de $J$ à partir de celui de $J_o$. Soit $G_o=]a,b[$ le plus grand intervalle ouvert contenu dans $\bigcup\limits_{w\in W(J_o)} (w+4^{-2N}E_u)$. L'extrémité droite de $G_o$, $b$, appartient à un $w+4^{-2N}E_u$ $(w\in W(J)\ba W(J_o))$ ; ce $w+4^{-2N}E_u$ n'a aucun point intérieur dans $G_o$ puisque les translatés de $4^{-2N}E_u$ par les éléments de $4^{-2N}V_{2N}$ sont disjoints presque partout, donc, par le lemme~2, il n'a aucun point dans $G_o$, donc $b$ est son extrémité gauche ; notons $w=w(b)$. Soit $G_o' =]a,b'[$ le plus grand intervalle ouvert contenant $G_o$ et contenu dans $G_o \cap (w(b)) +4^{-2N}
E_u)$, et définissons $w(b')$ à partir de $G_o'$ comme $w(b)$ à partir de $G_o$. D'après notre observation initiale,
$$
b'-b =w(b') -w(b) \ge 4^{-2N}[J|
$$
En poursuivant de la sorte, on définit $G_o'' = ]a,b''[$ et $w(b'')$, et
$$
b'' -b' \ge 4^{-2n}(J)
$$ 
et ainsi de suite. On construit de cette façon une chaîne de translatés de $4^{-2n}E_u$ qui va rejoindre $J_1$, puis comprendre les $w+ 4^{-2n}E_u$ $(w\in W(J_1))$. La chaîne se poursuit en se répétant au delà de $J_1$ jusqu'à la frontière de $E_u$. En partant de $J_1$ vers la gauche au lieu de $J_o$ vers la droite, la chaîne comprend les $w+4^{-2n}E_u$ $(w\in W(J_o))$ et se prolonge à gauche au delà de $J_o$ jusqu'à~0. La construction est périodique, de période $v_1-v_0$, et elle réalise un pavage du plus grand intervalle $]0,x[$ contenu dans $E_u$, périodique et de période $v_1-v_o$. Ce pavage se prolonge à droite en remplaçant $E_u$ par $4E_u$, $4^2E_u$ etc\dots\ Ainsi $\R^+$ est recouvert par un pavage périodique, de période $v_1-v_o$ dont les éléments s'écrivent $w+4^{-2n}E_u$~avec
$$
w= 4^{-2n} \sum_{k\ge 0} (\varepsilon_k+u\varepsilon_k')4^k\,,
$$
la somme étant finie et les $\varepsilon_k$ et $\varepsilon_k'$ égaux à 0 ou 1. Ces $w$ se répartissent en un nombre fini de progressions arithmétiques de raison $v_1-v_o$. Au moins une de ces progressions contient deux $\sum \varepsilon_k 4^k$ différents, dont la différence est donc un multiple de $v_1-v_o$. De même il existe deux $u \sum \varepsilon_k' 4^k$ dont la différence est un multiple de $v_1-v_o$. Il s'ensuit que $u$ est rationnel. 

Cela achève la démonstration du théorème~1.

\section{Interprétation du théorème 1 sur le modèle réduit, et complément}

Le modèle réduit est l'ensemble $B$ défini par (1), (2), (3), et ses sections horizontales sont les $B_h$ de (4). Puisque $B_o=E$ et $B_1=E'$, bornons--nous à $0<h<1$. Rappelons (5),~que
$$
\frac{u}{2} = \frac{h}{1-h}\,,\ \ h=\frac{u}{u+2}\,.
$$
Dire que $u$ est irrationnel, c'est dire que $h$ est irrationnel ; alors $mes_1 B_h=0$. Si $h=\frac{p}{q}$ irréductible, $u=\frac{2p}{q-p}$. Dire que $u$ est le rapport de deux nombres impairs, c'est dire que $p$ et $\frac{q-p}{2}$ sont impairs ; dans ce cas, $mes_1 B_h=0$ et $\dim B_h<1$. Il en est de même si $p^*$ et $\big(\frac{q-p}{2}\big)^*$ sont impaires. Si $p^*$ ou $\big(\frac{q-p}{2}\big)^*$ est paire (il n'est pas possible que les deux le soient), $B_h$ contient un intervalle et c'est l'adhérence de son intérieur. Ces propriétés sont invariantes par le changement de $u$ en $4u$, donc par le changement de $h$ en~$\frac{4h}{3h+1}$.

Nous avons vu de deux manières que $B$ est un ensemble de Besicovitch (au sens large). Pour obtenir un ensemble de Besicovitch au sens strict (le sens usuel), on peut faire pivoter $B$ de $k\frac{\pi}{4}$ autour de 0 $(k=1,2,3,4,5,6,7,8)$ et prendre la réunion des huit ensembles ainsi obtenus.

On peut aussi écraser $B$ sur la droite réelle et prendre la réunion des $B$ écrasés, complétée par le segment réel $[0,\frac{2}{3}]$. Précisons : il s'agit de l'adhérence~de
$$
\bigcup_{n=1}^\infty A_n(B)\,,
$$
où $A_n$ est la transformation affine $(x,y) \rg (x,\alpha_ny)$, $(\alpha_n)$ étant une suite positive quelconque tendant en décroissant vers~0.

La démonstration du théorème 1 contient des informations intéressantes sur la mesure $\mu$ dont la transformée de Fourier, écrite en (17), est
$$
\hat{\mu}(t) = \prod_{j=1}^\infty \cos 4^{-j} \frac{t}{2}\cos 4^{-j} \frac{ut}{2}
$$
et sur la réunion $V$ de tous les $V_n$, c'est--à--dire de l'ensemble
\begin{align}
V = \sum_{k\ge 0} (\varepsilon_k + u \varepsilon_k') 4^k
\end{align}
où la somme est finie et les $\varepsilon_k$ et $\varepsilon_k'$ valent 0 ou 1. Groupons--les.

\vskip2mm

\textsc{Théorème 2}.--- \textit{Si $u$ est irrationnel, $\mu$ est singulière et portée par un ensemble fermé de mesure nulle $(mes_1 E_u=0)$, et $V$ contient des points arbitrairement proches mais jamais confondus. Si $u$ est rationnel, $u=\frac{p}{q}$ irréductible, on a l'alternative suivante} :

$a)$ \textit{si $p^*$ et $q^*$ sont impaires, $\mu$ est singulière et portée par un ensemble fermé de dimension inférieure à $1$ $(\dim E_u<1)$, $\hat{\mu}(t)$ ne tend pas vers $0$ à l'infini, et $V$ contient des points multiples}

$b)$ \textit{si $p^*+q^*$ est impaire, c'est--à--dire l'un des deux pair et l'autre impair, $\mu$ est absolument continue, sa dérivée est bornée et son support est l'adhérence de son intérieur,  $V$ n'a pas de point multiple et c'est une réunion finie de progressions arithmétiques de même raison} :
$$
V=\alpha N + A\,,\ \alpha>0\,,\ \sharp A< \infty\,.
$$

On connait le comportement asymptotique de $\hat{\mu}(t)$ quand $u$ est rationnel : $\overline{\lim\limits_{t\rg \infty}} \hat{\mu}(t) >0$ dans le cas $a)$, et $\lim\limits_{t \rg \infty} \hat{\mu}(t)=0$ dans le cas $b)$. Il est possible que $\hat{\mu}(t)$ tende vers $0$ à l'infini quand $n$ est irrationnel, mais je ne sais pas le démontrer.

\section{Variantes}

Commençons par des variantes très proches du modèle réduit, pour aller ensuite vers des variantes hors du plan.

\subsection{} Au lieu de (1) et (2) on peut partir de
\begin{align}
& E= \displaystyle \Big\{ \sum_{j=1}^\infty \eta_j 9^{-j}\,,\ \eta_j \in \{0,1,2\}\Big\}\nonumber\\
&E' = 3E+i\nonumber
\end{align}
$B$ étant toujours défini comme la réunion des segments de droites joignant un point de $E$ et un point de $E'$ (formule~(3)).

L'étude est la même. On voit que $E+3E=[0,1]$, et~que
$$
3E-E\! =\! \Big\{\sum_{j=1}^\infty (3\eta_j' - \eta_j) 9^{-j}\Big\} \!=\! 
\Big\{\sum_{j=1}^\infty (3\eta_j' +(2-\eta_j)-2)9^{-j}\Big\}\! =\! \Big[-\frac{1}{4},\frac{3}{4}\Big]\,,
$$
donc les segments de droites constituant $B$ ont toutes les directions dans un certain angle. Les $E_u =E+uE$ sont de mesure nulle $(mes_1 E_u=0)$ pour tous les $u$ irrationnels, on le voit en transcrivant les lemmes~1, 2, 3. Reste le cas où $u$ est rationnel, $u=\frac{p}{q}$ irréductible. On définit $p^*$ et $q^*$ comme les premiers chiffres non nuls à partir de la droite dans les écritures de $p$ et de $q$ en base $9$, on transcrit de manière évidente les définitions des $V_n$, de $V$ et de la mesure $\mu$, et on a l'analogue des théorèmes~1 et 2 en remplaçant comme il convient les conditions sur $p^*$ et $q^*$, à savoir \og$3|q^*$ ou $3|p^*$\fg\ pour assurer que $V$ n'a pas de point multiple, et \og$3 \nmid q^*$ et $3\nmid p^*$\fg\ pour entraîner $\overline{\lim\limits_{t\rg \infty}} |\hat{\mu}(t)| >0$.

\subsection{} Plus généralement, quand $r$ est un nombre premier, on peut partir de
\begin{align}
& E= \displaystyle \Big\{ \sum_{j=1}^\infty \eta_j r^{-2j}\,,\ \eta_j \in \{0,1,\ldots r-1\}\Big\}\nonumber\\
&E' = rE+i\nonumber
\end{align}
et définir $E_u = E+ uE$, $V_n = \Big\{\sum\limits_{k=0}^{n-1} \eta_k r^{2k}\Big\}$ et $V=\cup V_n$, et $\mu$ comme mesure des probabilité canonique sur $E_u$,~avec
\begin{align}
&\hat{\mu}(t) = \displaystyle \prod_{j=1}^\infty C(r^{-2j}t)\,,\ C(t) =D(t) D(ut)\,, 
\nonumber\\
&D(t) = 1+ e^{it} +e^{2it} +\cdots e^{i(r-1)t} = \dfrac{1-e^{irt}}{1-e^{it}}\,.
\nonumber
\end{align}
De nouveau, les segments joignant $E$ et $E'$ constituent un ensemble de Besicovitch dont les sections horizontales de hauteur $h$ sont de la forme
\begin{align}
B_h = (1-h) E_u + ih\,,\ u=\frac{rh}{1-h}\,.
\end{align}
Pour $u=\frac{p}{q}$, $p^*$ et $q^*$ représentent les premiers chiffres non nuls en partant de la droite dans les écritures de $p$ et $q$ en base $r^2$. Les théorèmes~1 et 2 se transcrivent~ainsi.

\vskip2mm

\textsc{Théorème 3}.--- \textit{Si $u$ est irrationnel, $mes_1 E_u=0$, $\mu$ est singulière et $V$ contient des points arbitrairement proches mais jamais confondus. Si $u =\frac{p}{q}$ irréductible on a l'alternative suivante :}

$a)$ \textit{si $r$ ne divise ni $p^*$ ni $q^*$, $\dim E_u<1$, $\mu$ est singulière et $\overline{\lim\limits_{t\rg\infty}} |\hat{\mu}(t)|>0$, et $V$ contient des points multiples}

$b)$ \textit{si $r$ divise $p^*$ ou $q^*$, $E_u$ est l'adhérence de son intérieur, $\mu$ est absolument continue et sa dérivée est bornée, et $V$ est simple et de la forme}
$$
V = \alpha \N +A\,,\ \alpha>0\,,\ A\ \textrm{fini}\,.
$$

\subsection{} On peut aussi jouer avec deux nombres premiers différents, $r$ et $s$, en considérant
\begin{align}
E^r & = \displaystyle\Big\{\sum_{j=1}^\infty \eta_j(rs)^{-j}\,,\ \eta_j \in \{0,1,\ldots, r-1\}\Big\} \nonumber\\
E^s & = \displaystyle\Big\{\sum_{j=1}^\infty \eta_j'(rs)^{-j}\,,\ \eta_j \in \{0,1,\ldots, s-1\}\Big\} \nonumber\\
E_u &= E^r +u E^s \ \ (u>0)\nonumber\\
V & = \displaystyle\Big\{\sum_{k\ge0} \eta_k(rs)^{k}\,,\ \eta_k \in \{0,1,\ldots, r-1\}\Big\} \nonumber\\
V' & = \displaystyle\Big\{\sum_{k\ge0} \eta_k'(rs)^{k}\,,\ \eta_k \in \{0,1,\ldots, s-1\}\Big\}\,, \nonumber
\end{align}
les derniers $\sum$ représentant des sommes finies. Maintenant $\mu$, la mesure de probabilité canonique sur $E_u$, a pour transformée de Fourier
\begin{align}
\hat{\mu}(t) &= \prod_{j=1}^\infty C((rs)^{_j}t)\,,\ C(t) = D_r(t)D_s(ut)\,,\nonumber\\
D_r(t) &= 1+e^{it}+\cdots +e^{i(r-1)t}\nonumber\\
D_s(t) &= 1+e^{it}+\cdots +e^{i(s-1)t}\,.\nonumber
\end{align}

L'ensemble de Besicovitch associé est la réunion des segments joignant les ensembles $E^r$, porté par la droite $y=0$, et $rE^s+i$, porté par la droite $y=1$. On remarque~que
$$
E^r + r E^s = E_r = [0,1]\,,
$$
ce qui assure que les segments ont toutes les directions dans un ensemble ouvert de directions, les sections horizontales sont données par (19) en fonction des $E_u$, et les propriétés de $E_u$, $\mu$ et $V$ s'expriment par l'extension suivante du théorème~3.

\vskip2mm

\textsc{Théorème 4}.--- \textit{Même énoncé que le théorème $3$ pour $u$ irrationnel. Pour $u$ rationnel, $u=\frac{p}{q}$ irréductible, on désigne par $p^*$ $\resp q^*$ le premier chiffre non nul dans l'écriture de $p^*$ $\resp q^*$ en base $rs$. Même  énoncé que le théorème~$3$ en prenant comme hypothèse dans $a)$ que $r$ ne divise pas $p^*$ et que $s$ ne divise pas $q^*$, et dans $b)$ l'hypothèse opposée, c'est--à--dire que $r$ divise $p^*$ ou que $s$ divise~$q^*$.}

\subsection{}Si on laisse de coté les ensembles de Besicovitch, on peut étudier pour eux--mêmes les ensembles $E_u$ de la forme
$$
E_u = \Big\{\sum_{j=1}^\infty \eta_j r^{-j}\,,\ \eta_j \in \{0,1,\ldots,r-2,u\}\Big\} \quad (u>0)\,.
$$
L'exemple type est $r=3$, et c'est le cas traité dans l'article de Kenyon \cite{Keny}. Son résultat principal est le suivant

\vskip2mm

\textsc{Théorème de Kenyon} \cite{Keny} p. 224.--- \textit{Pour avoir $mes_1 E_u>0$, il faut et suffit que $u$ soit rationnel, égal à $\frac{p}{q}$ irréductible avec $p+q=0 \mod 3$. Si $u=\frac{p}{q}$ irréductible avec $p+q\not=0 \mod 3$, $\dim E_u<1$.}

\subsection{Le cas multidimensionnel}

On peut construire un ensemble de Besicovitch dans $\R^{d+1} $ en partant de $E^d$, $E$ étant défini en 4.2, ou à partir de $(E^r)^d$ et $(E^s)^d$, $E^r$ et $E^s$ étant définis en 4.3. Je me bornerai au cas $r=s=2$, c'est--à--dire à l'ensemble $E$ de départ, donné par (1), en considérant~$E^d$.

J'écrirai $\R^{d+1} =\R_1 \times \R_2\times\cdots \times \R_d \times \R_{d+1}$, et je désignerai par $i$ le vecteur unitaire de $\R_{d+1}$. L'ensemble $B$ et ses sections \og horizontales\fg\ $B_h$ seront encore définis par (3) et (4), avec ici $E^d$ à la place de $E$, et $(2E)^d+i$ à la place de $E'$. L'enveloppe convexe de $E^d$ est le cube $\big[0,\frac{1}{3}\big]^d$, et celle de $(2E)^d+i$ le cube $\big[0,\frac{2}{3}\big]^d+i$. Les segments de droite joignant un point de $E^d$ à un point de $(2E)^d+i$, dont la réunion est $B$, ont toutes les directions des segments joignant les deux cubes. Pour montrer que $B$ est un ensemble de Besicovitch, il suffit de montrer que $mes_d F_u=0$ pour presque tout $u$,~avec
\begin{align}
F_u &= E^d +u E^d \nonumber\\
&= \displaystyle\Big\{\sum_{j=1}^\infty (\varepsilon_j +u\varepsilon_j') 4^{-j}\,,\ \varepsilon_j \in \{0,1\}^d\,,\ \varepsilon_j' \in \{0,1\}^d\Big\}\nonumber\\
&= \displaystyle\Big\{\sum_{j=1}^\infty a_j  4^{-j}\,,\ a_j \in \{0,1,u,1+u\}^d\Big\}\,.\nonumber
\end{align}
En fait, on a $mes_d F_u=0$ pour tout $u$ irrationnel.

\vskip2mm

\textsc{Théorème 5} (Transcription du théorème 1).--- \textit{Si $u=\frac{p}{q}$ avec $p^*+q^*$ impair, $F_u$ est l'adhérence de son intérieur. Si $u=\frac{p}{q}$ avec $p^*$ et $q^*$ impairs, $mes_d F_u=0$ et $\dim F_u<d$. Si $u$ est irrationnel, $mes_d F_u=0$.}

\vskip2mm

La preuve est une retranscription de celle du théorème 1. Les lemmes~1 et 2 se retranscrivent verbatim. Le lemme~3 nécessite de remplacer au départ les intervalles par des boules, puis par des pavés. Je me borne à $d=2$ pour expliquer la situation. On considère $\R^2$ euclidien, on pose
$$
V_n = \Big\{ \sum_{k=0}^{n-1} b_k 4^k\,,\ b_k \in \{0,1,u,1+u\}^2\Big\}\,.
$$
On suppose que $F_u$ contient un disque $D$ de rayon $\rho$. La distance de deux points de $V_n$ est alors minorée par $\rho$. Pour $N$ assez grand, $D$ contient deux morceaux disjoints du recouvrement de $F_u$ par les $v+4^{-N} F_u(v\in 4^{-N}V_N)$ soit $v_o +4^{-N} F_u$ et $v_1 +4^{-N}F_u$, avec la condition additionnelle que $v_1$ soit sur la même horizontale que $v_o$ et à sa droite, ou sur la même verticale et au--dessus. Le premier, $v_o +4^{-N}F_u$, contient un rectangle ouvert $J_o$ dont un coté est parallèle à $v_1-v_o$, et le translaté de $J_o$ par $v_1-v_o$, $J_1$, occupe dans $v_1 +4^{-N}F_u$ la même position que $J_o$ dans $v_o +4^{-N}F_u$. Soit $B$ le coté de $J_o$ faisant face à $J_1$. Chaque point de $B$ appartient à un $w+4^{-2N}F_u(w\in 4^{-2N}V_{2N})$ disjoint de $J_o$, et c'est un point extrémal de ce $w+4^{-2N}V_{2N}$ dans la direction $v_o-v_1$ (parmi d'autres peut--être). On réunit ces $w+4^{-2N}V_{2N}$ à ceux qui recouvrent $J_o$ de façon minimale. L'ensemble obtenu contient un translaté de $J_o$ dans la direction $v_1-v_o$, $J'$, dont le coté qui fait face à $J_1$, $B'$, est le plus proche possible de $J_1$. En partant de $J'$ et $B'$ au lieu de $J_o$ et $B$, on définit $J''$ et $B''$ et ainsi de suite. Les segments $B$, $B'$, $B''$,\dots ne peuvent pas s'accumuler parce que les $w\in 4^{-2N}V_{2N}$ sont à des distances mutuelles $\ge 4^{-2N}\rho$. La suite $J_o$, $J'$, $J''$\dots se prolonge donc jusqu'à ce que le dernier terme coïncide avec $J_1$, et on peut la prolonger indéfiniment par périodicité de vecteur $v_1-v_o$. Le fait que $u$ est rationnel en dérive comme en dimension~1.

En partant de $J_1$ vers $J_o$ on a une chaîne de rectangles qui s'arrête au voisinage du bord de $\R^{+2}$, à savoir quand on se trouve à une distance de ce bord inférieure à $4^{-N}\rho$. En alternant les $v_1-v_o$ horizontaux et verticaux, on obtient un pavage doublement périodique de $[\varepsilon,\infty[^2$ par des $w+4^{-2N}F_u$, $w\in 4^{-2N}V$ avec $V = \bigcup\limits_n V_n$, pour~$N>N(\varepsilon)$.

L'étude pour $r$ rationnel se ramène au cas $d=1$.

\section{Commentaires sur les références}

L'histoire des ensembles de Besicovitch et leurs liens avec l'analyse harmonique, l'arithmétique et la combinatoire se trouvent excellemment décrits dans \cite{Laba}. Besicovitch les avait introduits dès 1918 pour montrer que les intégrales de Riemann ne se prêtent pas à un théorème de Fubini. Il en a repris la description en 1928 dans \cite{BesSam1} comme réponse au problème de Kakeya, puis en 1938 dans \cite{BesSam2} à l'occasion d'un film montrant la construction. La fin de cet article énonce une propriété remarquable des ensembles de droites dans le plan : on associe à chaque droite $D$ son pôle $P$ par rapport à un cercle fixé, et ainsi à un ensemble de droites $D_\alpha$ $(\alpha\in A)$ un ensemble de points $P_\alpha(\alpha\in A)$, soit $E$. Cette propriété est redonnée, et démontrée, dans \cite{BesSam3} :

\vskip2mm

\textsc{Théorème} (\cite{BesSam3} p. 706).--- \textit{Si l'ensemble $E$ est irrégulier (purement non rectifiable), l'ensemble $\bigcup\limits_{\alpha\in A} D_\alpha$ est de mesure plane nulle. S'il est régulier (rectifiable), la réunion des $D_\alpha$ est de mesure plane infinie.}
\vskip2mm

Besicovitch ajoute dans \cite{BesSam2} que la construction d'un ensemble $E$ irrégulier tel que les $D_\alpha$ aient toutes les directions ne présente pas de difficulté. Cela occupe quelques lignes dans \cite{BesSam3}, et se trouve expliqué complètement dans~\cite{BesSam4}.

L'exposé moderne le plus accessible sur les théorèmes de projection et sur les ensembles de Besicovitch et leurs généralisations se trouve dans le chapitre~18 du livre de Mattila \cite{Matt}. L'existence d'ensembles dans $\R^n$ constitués de plans ayant toutes les directions possibles, et de mesure nulle en dimension $n$, est un problème ouvert (\cite{Matt} p.~263). Mattila donne un exemple d'ensemble de Besicovitch par une variante de la méthode de dualité de Besicovitch 1964 \cite{BesSam4}, en associant au point $(a,b)$ la droite~$y=ax+b$.

La méthode de Kenyon exploitée ici dans la partie 2, et son théorème énoncé en 4.4 sont tirés de \cite{Keny}. Kenyon part du tamis défini dans $\C$~par
$$
S = \Big\{ \sum_{j=1}^\infty s_j 3^{-j}\,,\ s_j \in \{0,1,i\}\Big\}
$$
qui est un ensemble de Sierpinski de dimension~1, et purement non--rectifiable. La projection qui amène $x+iy$ sur $x+uy$ amène $S$~sur
$$
S_u = \Big\{ \sum_{j=1}^\infty a_j 3^{-j}\,,\ a_j \in \{0,1,u\}\Big\}\,,
$$
l'ensemble considéré en 4.4. J'ai juste transcrit sa méthode pour l'étude des ensembles $E_u$. Le théorème~4.4 est un abrégé très incomplet des résultats de Kenyon sur~les~$S_u$.

L'article de Kahane et Salem de 1958 \cite{KahSal} traite dans sa dernière partie des mesures $\mu$ dont les transformée de Fourier sont les produits $\gamma_\xi(t)\gamma_\xi(\lambda t)$,~où
$$
\gamma_\xi(t) = \prod_o^\infty \cos \pi \xi^kt\,.
$$
Il considère les cas $\xi>\frac{1}{4}$ ($\mu$ absolument continue pour presque tout $\lambda$) et $\xi<\frac{1}{4}$ ($\mu$ singulière pour tout $\lambda$), et laisse de coté les cas $\xi=\frac{1}{4}$. C'est le cas traité ici par le théorème 2.

Mon article de 1969 \cite{Kah} nécessitait une correction. Elle est faite.

\eject

\vskip4mm

\hfill\begin{minipage}{6,5cm}
Jean--Pierre Kahane

Laboratoire de Math\'ematique

Universit\'e Paris--Sud, B\^at. 425

91405 Orsay Cedex

13 juin 2012

\textsf{Jean-Pierre.Kahane@math.u-psud.fr}

\end{minipage}

\end{document}